 
\magnification=\magstep1       
\hsize=5.9truein                     
\vsize=8.5truein                       
\parindent 0pt
\parskip=\smallskipamount
\mathsurround=1pt
\hoffset=.25truein                     
\voffset=2\baselineskip               
%
%
\def\today{\ifcase\month\or
  January\or February\or March\or April\or May\or June\or
  July\or August\or September\or October\or November\or December\fi
  \space\number\day, \number\year}
%
 at 10truept
%
\newcount\dispno      
\dispno=1\relax       
\newcount\refno       
\refno=1\relax        
\newcount\citations   
\citations=0\relax    
\newcount\sectno      
\sectno=0\relax       
\newbox\boxscratch    
%

%
%
%
\def\Section#1#2{\global\advance\sectno by 1\relax%
\label{Section\noexpand~\the\sectno}{#2}%
\smallskip
\goodbreak
\setbox\boxscratch=\hbox{\bf Section \the\sectno.~}%
{\hangindent=\wd\boxscratch\hangafter=1
\noindent{\bf Section \the\sectno.~#1}\nobreak\smallskip\nobreak}}
%
\def\sqr#1#2{{\vcenter{\vbox{\hrule height.#2pt
              \hbox{\vrule width.#2pt height#1pt \kern#1pt
              \vrule width.#2pt}
              \hrule height.#2pt}}}}
\def\square{$\mathchoice\sqr34\sqr34\sqr{2.1}3\sqr{1.5}3$}
\def\endproof{~~\hfill\square\par\medbreak}
\def\noproof{~~\hfill\square}
%
%
\def\proc#1#2#3{{\hbox{${#3 \subseteq} \kern -#1cm _{#2 /}\hskip 0.05cm $}}}
\def\propcont{\mathchoice\proc{0.17}{\scriptscriptstyle}{}
                         \proc{0.17}{\scriptscriptstyle}{}
                         \proc{0.15}{\scriptscriptstyle}{\scriptstyle }
                         \proc{0.13}{\scriptscriptstyle}{\scriptscriptstyle}}
%

%
\def\normalin{\hbox{\raise0.045cm \hbox
                   {$\underline{\triangleleft }$}\hskip0.02cm}}
%
%
\def\'#1{\ifx#1i{\accent"13 \i}\else{\accent"13 #1}\fi}
%
%
%
\def\semidirect{\rlap{$\times$}\kern+7.2778pt \vrule height4.96333pt
width.5pt depth0pt\relax\;}
%
%
\def\prop#1#2{\noindent{\bf Proposition~\the\sectno.\the\dispno. }%
\label{Proposition\noexpand~\the\sectno.\the\dispno}{#1}\global\advance\dispno 
by 1{\it #2}\smallbreak}
\def\thm#1#2{\noindent{\bf Theorem~\the\sectno.\the\dispno. }%
\label{Theorem\noexpand~\the\sectno.\the\dispno}{#1}\global\advance\dispno
by 1{\it #2}\smallbreak}
\def\cor#1#2{\noindent{\bf Corollary~\the\sectno.\the\dispno. }%
\label{Corollary\noexpand~\the\sectno.\the\dispno}{#1}\global\advance\dispno by
1{\it #2}\smallbreak}
\def\defn{\noindent{\bf
Definition~\the\sectno.\the\dispno. }\global\advance\dispno by 1\relax}
\def\lemma#1#2{\noindent{\bf Lemma~\the\sectno.\the\dispno. }%
\label{Lemma\noexpand~\the\sectno.\the\dispno}{#1}\global\advance\dispno by
1{\it #2}\smallbreak}
\def\rmrk#1{\noindent{\bf Remark~\the\sectno.\the\dispno.}%
\label{Remark\noexpand~\the\sectno.\the\dispno}{#1}\global\advance\dispno
by 1\relax}
\def\proof{\noindent{\it Proof: }}
\def\numbeq#1{\the\sectno.\the\dispno\label{\the\sectno.\the\dispno}{#1}%
\global\advance\dispno by 1\relax}

\def\comm#1,#2{\left[#1{,}#2\right]}
\newdimen\boxitsep \boxitsep=0 true pt
\newdimen\boxith \boxith=.4 true pt 
\newdimen\boxitv \boxitv=.4 true pt
\gdef\boxit#1{\vbox{\hrule height\boxith
                    \hbox{\vrule width\boxitv\kern\boxitsep
                          \vbox{\kern\boxitsep#1\kern\boxitsep}%
                          \kern\boxitsep\vrule width\boxitv}
                    \hrule height\boxith}}
\def\square{\ \hbox{\vrule height7.5pt depth1.5pt width 6pt}\par}
\outer\def\square{\ifmmode\else\hfill\fi
   \setbox0=\hbox{} \wd0=6pt \ht0=7.5pt \dp0=1.5pt
   \raise-1.5pt\hbox{\boxit{\box0}\par}
}

\def\frac#1/#2{\leavevmode\kern.1em
              \raise.5ex\hbox{\the\scriptfont0 #1}\kern-.1em
              /\kern\.15em\lower.25ex\hbox{\the\scriptfont0 #2}}
\def\incnoteq{\lower.1ex \hbox{\rlap{\raise 1ex
     \hbox{$\scriptscriptstyle\subset$}}{$\scriptscriptstyle\not=$}}}
%
%


\def\propcontup{\bigcup\!\!\!\rlap{\kern+.2pt$\backslash$}\,\kern+1pt\vert}
%
%
%
\def\label#1#2{\immediate\write\aux%
{\noexpand\def\expandafter\noexpand\csname#2\endcsname{#1}}}
%
\def\ifundefined#1{\expandafter\ifx\csname#1\endcsname\relax}
%
%
\def\ref#1{%
\ifundefined{#1}\message{! No ref. to #1;}%
 \else\csname #1\endcsname\fi}
%
%
\def\refer#1{%
\the\refno\label{\the\refno}{#1}%
\global\advance\refno by 1\relax}
%
%
\def\cite#1{%
\expandafter\gdef\csname x#1\endcsname{1}%
\global\advance\citations by 1\relax
\ifundefined{#1}\message{! No ref. to #1;}%
\else\csname #1\endcsname\fi}
%
%
\font\bb=msbm10 
\font\bp=msbm10 at 8truept      
%
%
%

\def\Q{\hbox{\bb Q}}

\def\Z{\hbox{\bb Z}}                     \def\ZZ{\hbox{\bp Z}}

\def\Z{\hbox{\bb Z}}                     \def\ZZ{\hbox{\bp Z}}

\newread\aux
\immediate\openin\aux=\jobname.aux
\ifeof\aux \message{! No file \jobname.aux;}
\else \input \jobname.aux \immediate\closein\aux \fi
\newwrite\aux
\immediate\openout\aux=\jobname.aux
 
\font\smallheadfont=cmr8 at 8truept

\headline={\ifnum\pageno<2{\hfill}\else{\ifodd\pageno\rightheadline
\else\leftheadline\fi}\fi}
\def\leftheadline{\smallheadfont A. Magidin\hfil}
\def\rightheadline{\hfil\smallheadfont Nonsurjective epimorphisms in
decomposable varieties of groups}
 
\centerline{\bf Nonsurjective epimorphisms in decomposable varieties of groups}
\centerline{Arturo Magidin\footnote*{The author was
supported in part by a fellowship from the Programa de Formaci\'on y
Superaci\'on del Personal Acad\'emico de la UNAM, administered by the
DGAPA.}}
\centerline{\today}
\smallskip
{\parindent=20pt
\narrower\narrower
\noindent{\smallheadfont{Abstract. A full characterization of when a
subgroup~$\scriptstyle H$ of a group~$\scriptstyle G$ in a varietal
product $\scriptstyle {\cal NQ}$ is epimorphically embedded
in~$\scriptstyle G$ (in the
variety $\scriptstyle {\cal NQ}$) is given. From this,
a result of S.~McKay is derived , which states that if $\scriptstyle{\cal
NQ}$ has instances of nonsurjective epimorphisms, then
$\scriptstyle{\cal N}$ also has instances of nonsurjective
epimorphisms. Two partial converses to McKay's result are also given:
when~$\scriptstyle G$ is a finite nonabelian simple group; and
when~$\scriptstyle G$ is finite and $\scriptstyle {\cal Q}$ is a
product of varieties of nilpotent groups, each of which contains the
infinite cyclic~group.\par}}}
\bigskip
\medskip

\footnote{}{\noindent\smallheadfont Mathematics Subject
Classification:
08B25, 20E10, 20J99 (primary)}
\footnote{}{\noindent\smallheadfont Keywords:dominion, epimorphism,
varieties of groups, decomposable variety}

\Section{Introduction and notation}{intro}

Given a category~${\cal C}$, a map $f\colon G\to K$ in~${\cal C}$ is
an {\it epimorphism} if and only if it is right
cancellable. When~${\cal C}$ is a full subcategory of the category of
all algebras (in the sense of Universal Algebra) of a given type, it
is not hard to verify that if~$f$ is surjective, then it is an
epimorphism. The converse, however, does not necessarily hold. For
example, the embedding $\Z\hookrightarrow\Q$ is an epimorphism in the
category of rings, but it is not~surjective.

On the other hand, it is known that in the category of all groups,
epimorphisms are surjective; for an elementary proof of this fact we
direct the reader to~{\bf [\cite{episingroups}]}. Peter Neumann
proved~{\bf [\cite{pneumann}]} that in a full subcategory of~${\cal
G}roup$ in which all objects are solvable groups, and which is closed
under taking quotients and subgroups, all epimorphisms are
surjective. Susan McKay {\bf [\cite{mckay}]} later extended this
result. On the other hand, an example of B.H.~Neumann that appears
in~${\bf [\cite{pneumann}]}$ shows that there are varieties of groups
where there are nonsurjective epimorphisms. Specifically, the
embedding $A_4\hookrightarrow A_5$ is an epimorphism in~${\rm
Var}(A_5)$. For other examples of nonsurjective epimorphisms in
varieties of groups, we direct the reader to~{\bf
[\cite{simpleprelim}]} and~{\bf [\cite{mythesis}]}.

In this work we study the question of when a varietal product ${\cal
NQ}$, where ${\cal N}$ and~${\cal Q}$ are varieties of groups, has
instances of nonsurjective~epimorphisms.

Isbell introduced the concept of dominions in~${\bf
[\cite{isbellone}]}$ to study epimorphisms. Given a category~${\cal
C}$ as above, and an algebra $A\in {\cal C}$, Isbell defines for a
subalgebra~$B$ of~$A$ the {\it dominion of~$B$ in~$A$ (in the
category~${\cal C}$)} as the intersection of all equalizer subalgebras
of~$A$ containing~$B$. Explicitly,
$${\rm dom}_A^{\cal C}(B)=\Bigl\{a\in A\bigm| \forall C\in {\cal
C},\;\forall f,g\colon A\to C,\ {\rm if}\ f|_B=g|_B{\rm\ then\ }
f(a)=g(a)\Bigr\}.$$

Note that given an arbitrary morphism of algebras $\varphi\colon A'\to
A$, $\varphi$ is an epimorphism in~${\cal C}$ if and only if ${\rm
dom}_A^{\cal C}(\varphi(A'))=A$. Also, in a variety of algebras (and
in any full category of algebras which is closed under subalgebras and
quotients) we can factor $\varphi$ into a
surjection $A'\mapsto \varphi(A')$ and an embedding $\varphi(A')
\hookrightarrow A$. The description of surjective maps by means of
congruence relations is well developed, so we may reduce the study of
epimorphisms (at least in certain categories) to the study of~dominions.

For the basic properties of dominions in varieties of groups we refer
the reader to~{\bf [\cite{domsmetabprelim}]}. We recall the most
important properties: ${\rm dom}_G^{\cal V}(-)$ is a closure operator
on the lattice of subgroups of~$G$; the dominion construction respects
finite direct products; and the dominion construction respects
quotients. That is, if~$H$ is a subgroup of~$G\in{\cal V}$, where
${\cal V}$ is a variety, and
$N\triangleleft G$, with $N$ contained in~$H$, then
$${\rm dom}_{G/N}^{\cal V}(H/N) = {\rm dom}_G^{\cal V}(H)\bigm/ N.$$

Groups will be written multiplicatively unless otherwise stated. Given
a group~$G$, the identity element of~$G$ will be written $e_G$,
although we will omit the subscript if it is understood from
context. All maps will be assumed to be group morphisms unless
otherwise~specified. Given a group~$G$ and a subgroup~$H$, $N_G(H)$
denotes the normalizer of~$H$ in~$G$; that is, the subgroup of all elements
$g\in G$ such that $g^{-1}Hg=H$.

Given two groups, $A$ and~$B$, we write $A\wr B$ to denote the regular
wreath product of~$A$ and~$B$; this is the semidirect product of
$|B|$-copies of~$A$ (indexed by the elements of~$B$) by~$B$, with~$B$
acting on the index set via the regular right action. The elements
of~$A\wr B$ are written as $b\phi$, where $b\in B$ and $\phi\colon
B\to A$ is a set-theoretic function (that is, an element of~$A^{B}$). 

Recall that a variety of groups is a full subcategory of~${\cal
G}roup$ which is closed under taking quotients, subgroups, and
arbitrary direct products. For the basic properties and facts about
varieties of groups, we direct the reader to Hanna Neumann's excellent
book~{\bf
[\cite{hneumann}]}. We will denote the variety of all groups by~${\cal
G}$, and the variety consisting only of the trivial group by~${\cal
E}$.

Given two varieties of groups ${\cal N}$
and~${\cal Q}$, their product variety ${\cal NQ}$ is the variety of
all groups which are an extension of an ${\cal N}$-group by a ${\cal
Q}$-group; that is, all groups~$G$ with a normal subgroup
$N\triangleleft G$ such that $N\in {\cal N}$ and $G/N \in {\cal
Q}$. The product ${\cal NQ}$ is easily seen to be a variety, say by
using Birkhoff's HSP theorem~{\bf [\cite{birkhoff}]}. Multiplication
of varieties is~associative.

The semigroup of varieties of groups has the structure of a cancellation
semigroup with 0 and 1. The zero element is the variety of all groups,
while the identity is the trivial variety,
consisting only of the trivial group. 
We will say that a variety ${\cal V}$ is {\it nontrivial} iff ${\cal
V}\not={\cal G}$ and ${\cal V}\not={\cal E}$, and we will call it {\it
trivial} otherwise. A variety ${\cal V}$ {\it factors nontrivially}
(or is {\it decomposable}) if it can be expressed as the product of
two nontrivial~varieties.

Furthermore, every variety other than~${\cal G}$ can be uniquely
factored as a product of a finite number of indecomposable varieties
(in a unique order), with~${\cal E}$ having the empty factorization,
so that the semigroup with neutral element of varieties other
than~${\cal G}$ is freely generated by the indecomposable
varieties. See Theorems~21.72, 23.32 and~23.4 in {\bf
[\cite{hneumann}]}.

Given a variety ${\cal V}$ and a group~$G$ (not necessarily in~${\cal
V}$), we will denote by~${\cal V}(G)$ the verbal subgroup of~$G$
associated to~${\cal V}$. This is the fully invariant subgroup
generated by all values of the words ${\bf v}$ which are laws
of~${\cal V}$. In particular, $G\in {\cal V}$ if and only if ${\cal
V}(G)=\{e\}$. We also note the universal property associated to ${\cal
V}(G)$: for any normal subgroup $N\triangleleft G$, $G/N\in {\cal V}$
if and only if~${\cal V}(G)\subseteq N$.

In \ref{neumannresult} we will recall the result of Peter M.~Neumann
mentioned above. We will also recall a result from~{\bf
[\cite{prodvarprelim}]}, which gives upper and lower bounds for the
dominion of a subgroup in a decomposable~variety.  Next, in
\ref{charactofepis} we will characterize when a subgroup~$H$ of~$G$ is
epimorphically embedded in~$G$ in the variety~${\cal NQ}$ in terms of
the epimorphisms of~${\cal N}$, the laws of~${\cal Q}$, and the
internal structure of~$G$. \ref{generalnonsurjinprod}
is the main result of this work. Finally, in \ref{finitepis} and
\ref{partialtwo} we will prove two partial converses to
McKay's~theorem.

The contents of this work are part of the author's doctoral
dissertation, which was conducted under the direction of Prof.~George
M.~Bergman, at the University of California at~Berkeley. It is my very
great pleasure to record and express my deep gratitude and
indebtedness to Prof.~Bergman. His advice and suggestions were
invaluable, and improved this work in ways too numerous to mention
explicitly. He also helped correct many mistakes; any errors that
remain, however, are my own~responsibility.

\Section{Preliminary results}{neumannresult}

In this section we generalize a theorem of P.M. Neumann about
nonsurjective epimorphisms in certain classes of~groups. The proof
follows Neumann's proof very~closely.

\thm{generalpneumann}{(P.M.~Neumann {\bf [\cite{pneumann}]}) Let
${\cal X}$ be a full subcategory of ${\cal
G}roup$, where ${\cal O}b({\cal X})$ is closed under taking quotients
and subgroups of objects of~${\cal X}$. Let $G\in {\cal X}$, and let
$H$ be a proper subgroup of~$G$. Suppose that there is a normal
subgroup $N\triangleleft G$ such that $N$ is solvable and $NH=G$. Then
$${\rm dom}_G^{\cal X}(H)\propcont G.$$}

\proof Let $G\in {\cal X}$, and let $H$ be a proper subgroup
of~$G$. Let 
$$\{e\}=N_m\triangleleft N_{m-1}\triangleleft\cdots\triangleleft N_0=N$$
be a normal series for~$N$ such that $N_i\triangleleft G$ for each
$i$, and $N_{i}/N_{i+1}$ is abelian for $i=0,1,\ldots,m-1$ (we can
obtain such a series by letting $N_i$ be the $i$-th derived subgroup
of~$N$).

Since $N_mH=H\propcont G$, and $N_0H=NH=G$, there exists $k$ such that
$N_{k+1}H\propcont G$ and $N_{k}H=G$. It will 
suffice to show that ${\rm dom}_G^{\cal X}(N_{k+1}H)\propcont G$.

Let $H^*=N_{k+1}H$ and let $M=N_{k}\cap H^*$.  Then $M\triangleleft
H^*$, since $N_{k}\triangleleft G$; and $M\triangleleft N_{k}$,
because $N_{k+1}\subseteq M$, and so $M$ corresponds to
subgroup of $N_{k}/N_{k+1}$, which is abelian. Therefore,
$M\triangleleft N_{k}H^* = G$.

By assumption, $G/M\in {\cal X}$. Let $\theta\colon G\to G/M$ be the
quotient map.  Then $G/M$ is a semidirect product of~$\theta(N_k)$ by
$\theta(H^*)$, because $\theta(N_k)\triangleleft \theta(G)=G/M$, and the
kernel is contained in both $H^*$ and $N_k$, so 
$$\theta(H^*)\cap
\theta(N_k) = \theta(H^*\cap N_k) = \theta(M)=\{e\}.$$
Finally, note that $\theta(H^*)\theta(N_k) =
\theta(H^*N_k)=\theta(G)=G/M$.

If we now compare the map $\theta$ with the map obtained by composing
$\theta$ with the idempotent endomorphism of~$G/M$ with kernel
$\theta(N_k)$ and image $\theta(H^*)$, we see that
${\rm dom}_G^{\cal X}(H^*) = H^*$, and therefore we conclude that
${\rm dom}_G^{\cal X}(H)\subseteq H^*\propcont G$,
as~claimed.\endproof

\cor{solvthenepissurje}{(P.M.~Neumann {\bf [\cite{pneumann}]})
Let ${\cal X}$ be a full subcategory of ${\cal G}roup$, where all
objects in~${\cal X}$ are solvable groups, and ${\cal O}b({\cal X})$
is closed under taking quotients and subgroups of objects of~${\cal
X}$. Then all epimorphisms in~${\cal X}$ are~surjective.\noproof}

\cor{solvablevars}{If ${\cal V}$ is a variety of solvable groups, then
all epimorphisms are surjective in~${\cal V}$.\noproof}

\thm{mckaygeneralized}{(S. McKay {\bf [\cite{mckay}]}) Let ${\cal
V}={\cal NQ}$ be a variety, where ${\cal N}$ is a nontrivial
variety. Let $G\in {\cal V}$, and let $N\triangleleft G$, with
$N\in{\cal N}$ and $G/N\in {\cal Q}$. Then for all subgroups~$H$
of~$G$,
${\rm dom}_G^{\cal V}(H)\subseteq NH$. In particular,
$${\rm dom}_G^{\cal V}(H)\subseteq {\cal Q}(G)H.\eqno\noproof$$}

Finally, we recall the following result:

\thm{bigone}{(Theorem~3.12, {\bf [\cite{prodvarprelim}]}) Let ${\cal V}=
{\cal NQ}$ be a nontrivial factorization of~${\cal V}$, and let $G\in
{\cal V}$. Let $H$ be a subgroup of~$G$. If $D={\rm
dom}_{{\cal Q}(G)}^{\cal N}({\cal Q}(G)\cap H)$, then
$$\langle H,D\rangle = HD \subseteq {\rm dom}_G^{\cal NQ}(H).$$
Furthermore, if $N_G(D){\cal Q}(G)=G$, then
$${\rm dom}_G^{\cal NQ}(H)=HD$$
and ${\rm dom}_G^{\cal NQ}(H)\cap {\cal Q}(G)=D.$\noproof}

Recall that two varieties ${\cal N}$ and~${\cal Q}$ are {\it
disjoint} if and only if ${\cal N}\cap {\cal Q}={\cal E}$.

\cor{disjointvar}{Let ${\cal V}={\cal NQ}$, where ${\cal N}$ and
${\cal Q}$ are disjoint nontrivial varieties of groups. Let $G\in
{\cal N}$, and let $H$ be a subgroup of~$G$. Then
$${\rm dom}_G^{{\cal NQ}}(H) = {\rm dom}_G^{\cal
N}(H).\eqno\noproof$$}

Recall that B.H.~Neumann proved that the embedding $A_4\hookrightarrow
A_5$ is an epimorphism in the variety ${\rm Var}(A_5)$ (see Example~A
in~{\bf [\cite{pneumann}]}); since there are uncountably many
varieties of groups disjoint from~${\rm Var}(A_5)$, it follows from
\ref{disjointvar} that there are uncountably many varieties of groups
which contain instances of nonsurjective epimorphisms (see
Theorem~4.27 in~{\bf [\cite{prodvarprelim}]}).

\Section{Nonsurjective epimorphisms in decomposable varieties}{charactofepis}

What we can say about nonsurjective epimorphisms in ${\cal NQ}$ if we
drop the requirement that ${\cal N}\cap {\cal Q}={\cal E}$?

\thm{nonsurjinprod}{Let ${\cal V}={\cal NQ}$ be a variety of groups,
with ${\cal N}$ and~${\cal Q}$ nontrivial. Let $G\in{\cal V}$,
and let $H$ be a subgroup of~$G$. Then ${\rm
dom}_G^{\cal NQ}(H) = G$ if and only if $H{\cal Q}(G)=G$ and ${\rm
dom}_{{\cal Q}(G)}^{\cal N}\bigl(H\cap {\cal Q}(G)\bigr)={\cal
Q}(G)$.}

\proof Write $D={\rm dom}_{{\cal Q}(G)}^{\cal N}\bigl(H\cap {\cal
Q}(G)\bigr)$. We first prove the ``if'' part.

In this case we have
$$\eqalign{G &= H{\cal Q}(G)\qquad\qquad \hbox{(by hypothesis)}\cr
&= HD\qquad\qquad\quad\>\> \hbox{(by hypothesis)}\cr
&\subseteq {\rm dom}_G^{\cal NQ}(H)\qquad\> \hbox{(by
\ref{bigone})}\cr
&\subseteq G.\cr}$$
Therefore, $G={\rm dom}_G^{\cal NQ}(H)$, as claimed.

For the converse, first note that ${\rm dom}_G^{\cal NQ}(H)\subseteq
H{\cal Q}(G)$, by \ref{mckaygeneralized}. Therefore, $H{\cal
Q}(G)=G$. Since~$H$ normalizes itself, and everything normalizes
${\cal Q}(G)$, $H$ normalizes $H\cap {\cal Q}(G)$. Since the dominion
construction respects automorphisms, it follows that $H$ must
normalize~$D$. We conclude that
$H\subseteq N_G(D)$, and hence it follows that $N_G(D){\cal Q}(G)=G$. Applying
\ref{bigone} we have $${\rm dom}_G^{\cal NQ}(H) = \langle H,D\rangle = HD.$$

We also note that ${\rm dom}_G^{\cal
NQ}(H)\cap {\cal Q}(G) = D$. Since ${\rm dom}_G^{\cal NQ}(H)=G$ by
hypothesis,  $D={\cal Q}(G)$, as~claimed.\endproof

Since $\psi\colon H\to G$ is an epimorphism in the variety
${\cal V}={\cal NQ}$ if and only if $${\rm dom}_G^{\cal V}(\psi(H)) = G,$$
\ref{nonsurjinprod} tells us for which groups $G\in{\cal N}$ there is
a nonsurjective epimorphism with codomain~$G$, 
in terms of the varieties ${\cal N}$ and~${\cal Q}$, and the structure of
the group~$G$. We also get the following corollary:

\cor{nsinprodonlyiffirst}{(S.~McKay {\bf [\cite{mckay}]}) Let ${\cal
V}={\cal NQ}$ be a variety with
${\cal N}$ and ${\cal Q}$ nontrivial. If ${\cal V}$ has instances of
nonsurjective epimorphisms then ${\cal N}$ has instances of
nonsurjective epimorphisms.}

\proof By \ref{nonsurjinprod}, if a subgroup $H$ of~$G\in {\cal
V}$ is epimorphically embedded in~$G$ then $H\cap {\cal Q}(G)$ is
epimorphically embedded (in ${\cal N}$) in~${\cal Q}(G)$.  If $H$
contains ${\cal Q}(G)$, then $H=H{\cal Q}(G)=G$. Therefore, if~$H$
is a proper subgroup of~$G$, we must have $H\cap {\cal Q}(G)\propcont
{\cal Q}(G)$, so the 
embedding $H\cap {\cal Q}(G) \hookrightarrow {\cal Q}(G)$ is a
nonsurjective epimorphism in the variety~${\cal N}$.\endproof

In a way, \ref{nonsurjinprod} tells us that the existence of
nonsurjective epimorphisms is determined by the indecomposable
varieties. For given a variety ${\cal V}$ and a group $G\in{\cal V}$,
in order to find out if a subgroup $H$ of~$G$ is epimorphically embedded
in~$G$ we only need to do the following: factor ${\cal V}$ into a
product of indecomposable varieties ${\cal V}_1{\cal V}_2\cdots{\cal
V}_n$, and test to see whether $H{\cal V}_2\cdots{\cal V}_n(G)=G$ and
$H\cap {\cal V}_2\cdots{\cal V}_n(G)$ is epimorphically embedded in
${\cal V}_2\cdots{\cal V}_n(G)$ in the variety ${\cal V}_1$. 

On the other hand, carrying out these calculations may not be a trivial
matter. Specifically, Kle\u{\i}man has established that there does not
exist an algorithm that determines whether an arbitrary finitely based
variety decomposes into a product of two varieties (see Theorem~4.3
and introductory comments in~{\bf [\cite{kleimanbig}]}). 

Before continuing, we give a variant of \ref{nonsurjinprod}.

\lemma{epiinprodmoregen}{Let ${\cal V}$ be a variety of groups,
$G\in {\cal V}$ and $H$ a proper subgroup of~$G$ such that ${\rm
dom}_G^{\cal V}(H)=G$. Let $N\triangleleft G$ and suppose that ${\cal
A}$, ${\cal B}$ are two varieties such that $N\in {\cal A}$, $G/N\in
{\cal B}$, and ${\cal AB}\subseteq {\cal V}$. Then $HN=G$ and ${\rm
dom}_N^{\cal A}(H\cap N)=N$. In particular, ${\cal A}$ has instances
of nonsurjective epimorphisms.}

\proof Since $G\in {\cal AB}\subseteq{\cal V}$, it follows that ${\rm
dom}_G^{\cal AB}(H)=G$. By \ref{nonsurjinprod}, we must have $H{\cal
B}(G)= G$ and ${\rm dom}_{{\cal B}(G)}^{\cal A}(H\cap {\cal
B}(G))={\cal B}(G)$.

Since $G/N\in {\cal B}$, we also have ${\cal B}(G)\subseteq
N$. Hence
$$\eqalign{G&= H{\cal B}(G)\cr
&\subseteq HN\qquad\qquad\hbox{(since ${\cal B}(G)\subseteq N$)}\cr
&\subseteq G.\cr}$$
Therefore, $HN=G$, as claimed.

Since ${\cal B}(G)\subseteq N$, and ${\rm dom}_{{\cal B}(G)}^{\cal
A}(H\cap {\cal B}(G))={\cal B}(G)$, it follows that
$$\langle H\cap N,{\cal B}(G)\rangle =
(H\cap N){\cal B}(G) \subseteq {\rm dom}_N^{{\cal A}}(H\cap N).$$

We claim that in fact $N=(H\cap N){\cal B}(G)$, which will prove the
lemma. Indeed, let $n\in N$. Since $H{\cal B}(G)=G$, there exist $h\in
H$ and $b\in {\cal B}(G)$ such that $hb=n$. Therefore,
$h=nb^{-1}$. Since ${\cal B}(G)\subseteq N$, it follows that $h\in
N$. Therefore, $n\in (H\cap N){\cal B}(G)$, as~claimed.\endproof

From \ref{nonsurjinprod} and \ref{epiinprodmoregen}, we conclude our
main result:

\thm{generalnonsurjinprod}{Let ${\cal V}={\cal NQ}$ be a variety with
${\cal N}$ and ${\cal Q}$ nontrivial, and let $G$ be a group in ${\cal
V}$. For a subgroup $H$ of~$G$, the following are equivalent:}
\nobreak
{\it\parindent=25pt
\item{(i)} ${\rm dom}_G^{\cal V}(H)=G$.
\item{(ii)} $H$ is epimorphically embedded in~$G$ in the variety ${\cal
V}$.
\item{(iii)} $H{\cal Q}(G)=G$ and ${\rm dom}_{{\cal Q}(G)}^{\cal
N}(H\cap {\cal Q}(G))={\cal Q}(G)$.
\item{(iv)} For every normal subgroup~$N$ of~$G$ such that $N\in {\cal
N}$ and $G/N\in {\cal Q}$, $HN=G$ and ${\rm dom}_N^{\cal N}(H\cap
N)=N$.\par}

\proof By definition of dominion, (i) and~(ii) are equivalent. The
equivalence of (i) and~(iii) follows from
\ref{nonsurjinprod}. Clearly, (iv) implies~(iii). Finally, (i) implies~(iv) by
\ref{epiinprodmoregen}, setting ${\cal A}={\cal
N}$, and ${\cal B}={\cal Q}$.\endproof

\rmrk{needallnotjustone} We quickly note that
in~\ref{generalnonsurjinprod}(iv), it is not enough to
consider a single normal subgroup~$N$. We defer an example for~now.

We may ask whether the converse of \ref{nsinprodonlyiffirst}
holds. That is, if a variety ${\cal N}$ has instances of nonsurjective
epimorphisms, and ${\cal Q}$ is a nontrivial variety, does ${\cal NQ}$
have instances of nonsurjective epimorphisms? We will partially answer
this question in the next sections.

\Section{Finite epimorphisms}{finitepis}

We will say that a variety~${\cal N}$ has instances of {\it finite}
nonsurjective epimorphisms if there exists a finite group $G\in{\cal
N}$ and a proper subgroup~$H$ of~$G$ with ${\rm dom}_G^{\cal
N}(H) = G$. 

In this section we will prove two partial converses to
\ref{nsinprodonlyiffirst}. First, we will show that
if ${\cal N}$ has instances of finite nonsurjective epimorphisms, then
so does ${\cal NA}$, where ${\cal A}$ is the variety of abelian
groups; then we will show that the same holds if ${\cal A}$ is
replaced by any product of varieties of nilpotent groups, each of
which contains the infinite cyclic group. We need some preliminary
results regarding the commutator of a standard wreath product~first.

Let $\Z=\langle x\rangle$ denote the infinite cyclic group, which we
will write~multiplicatively. 

\lemma{commofwr}{Let $G$ be a group, and let $\Z$ be the
infinite cyclic group. Then the commutator subgroup
of $G\wr\Z$ is equal to $G^{\ZZ}$. In fact, letting $x$ denote a
generator of $\Z$, every element of~$G^{\ZZ}$ has the form $[\psi,x]$
for some $\psi\in G^{\ZZ}$.}

\proof The regular wreath $G\wr\Z$ is a semidirect product of $G^{\ZZ}$
with $\Z$. If we quotient out by $G^{\ZZ}$ we obtain a group isomorphic
to $\Z$, hence abelian. Therefore, $[G\wr\Z,G\wr\Z]\subseteq G^{\ZZ}$.

Let $\phi\colon \Z\to G$ be an arbitrary element of $G^{\ZZ}$. We claim
that there exists an element $\psi\in G^{\ZZ}$ such that $[\psi,x] =
\phi$ in~$G\wr\Z$.

We note that by definition of the wreath product, 
$$\psi^x(x^n) =
\psi(x^nx^{-1}) = \psi(x^{n-1}).$$
In particular, we have
$$[\psi,x](x^n) = \psi(x^n)^{-1}\psi(x^{n-1}).\eqno(\numbeq{basicident})$$

We define $\psi$ recursively. Choose $\psi(e_{\ZZ})$ arbitrarily. Assuming
we have defined $\psi$ at $x^n$, $n\geq 0$, since we want
$[\psi,x]=\phi$,  from (\ref{basicident}) we conclude that we must
define $$\psi(x^{n+1}) = \psi(x^n)\phi(x^{n+1})^{-1}.$$

Similarily, if we have defined $\psi$ at $x^{-n}$, for some $n\geq 0$,
again from $[\psi,x]=\phi$ and (\ref{basicident}), we must define
$\psi(x^{-n-1}) = \psi(x^{-n})\phi(x^{-n})$.

This defines $\psi$ recursively, and by construction $[\psi,x] =
\phi$, as desired. Therefore, $G^{\ZZ}$ is contained in $[G\wr\Z,G\wr\Z]$,
giving~equality.\endproof

\rmrk{neumannhadit} This result also follows from work of P.M.~Neumann
on the wreath product, specifically Corollary~5.3 in~{\bf
[\cite{neumannwreath}]}. 

\thm{commofwrn}{If $G$ is a group, then $$\left[(G\wr
\Z)^I,(G\wr\Z)^I\right] =
\bigl([G\wr\Z,G\wr\Z]\bigr)^I\eqno(\numbeq{eqofcomms})$$
for any set~$I$.}

\proof  Clearly, the left hand side is contained in the right hand
side of (\ref{eqofcomms}). Conversely, let $(\psi_i)_{i\in I}$ be an
element of $([G\wr\Z,G\wr\Z])^I$; thus $\psi_i\in [G\wr\Z,G\wr\Z] =
G^{\ZZ}$ for each $i\in I$. 

From \ref{commofwr} we know that every element of $[G\wr\Z,G\wr\Z]$ is
equal to the commutator of two elements; therefore, for each $i\in I$
there exist elements $x_i$, $y_i$ in $G\wr\Z$ such that $\psi_i =
[x_i,y_i]$. But this implies that
$$\Bigl[ \bigl(x_i\bigr)_{i\in I}, \bigl(y_i\bigr)_{i\in I}\Bigr] =
\Bigl(\bigl[x_i,y_i\bigr]\Bigr)_{i\in I} = \bigl(\psi_i\bigr)_{i\in
I}$$
which shows that the right hand side of~(\ref{eqofcomms}) is contained
in the left hand side, proving~equality.\endproof

\lemma{infforfinite}{Let ${\cal N}$ be a variety of groups, and let $G\in{\cal
N}$ be a finite group. If~$H$ is a subgroup of~$G$ and $I$ is an
arbitrary set, then $${\rm dom}_{G^I}^{\cal N}\left(H^I\right) =
\Bigl({\rm dom}_G^{\cal N}(H)\Bigr)^I.$$}

\proof Clearly, we have ${\rm dom}_{G^I}^{\cal N}(H^I)\subseteq
\bigl({\rm dom}_G^{\cal N}(H)\bigr)^I$.  To prove the reverse
inclusion, let $\phi\colon I\to {\rm dom}_G^{\cal N}(H)$ be an element
of $\bigl({\rm dom}_G^{\cal N}(H)\bigr)^I$.
 
For each $g\in G$, let
$$S_g=\{i\in I\mid \phi(i)=g\}.$$
 
Consider the group $K=G^{|{\rm dom}_G^{\cal N}(H)|}$, and let
$M=H^{|{\rm
dom}_G^{\cal N}(H)|}$. Since the number of direct factors is finite,
and dominions respect finite direct products, we have
$${\rm dom}_K^{\cal N}(M) = \bigr({\rm dom}_G^{\cal N}(H)\bigr)^{|{\rm
dom}_G^{\cal N}(H)|}.$$
 
Let $\eta\colon K\to G^I$ be the embedding that maps the
$g$-coordinate of~$K$ diagonally to the $S_g$ coordinate
of~$G^I$. That is, $\eta$ sends an element ${\bf z}\in K$ to the
element $\eta({\bf z})\in G^I$, where $\eta({\bf z})_i = {\bf z}_g$ if
and only if $i\in S_g$.
 
Clearly, $\eta(M)\subseteq H^I$, and therefore
$$\eta\bigl({\rm dom}_K^{\cal N}(M)\bigr) \subseteq {\rm
dom}_{G^I}^{\cal N}(H^I).$$
 
We claim that $\phi\in\eta({\rm dom}_K^{\cal N}(M))$. Indeed, let $\bf
z$ be the element of~$K$ given by ${\bf z}_g=g$ for each $g\in {\rm
dom}_G^{\cal N}(H)$. By construction of~$\eta$, $\eta({\bf z})=\phi$,
which proves the~claim.
 
Therefore, $\bigl({\rm dom}_G^{\cal N}(H)\bigr)^I\subseteq {\rm
dom}_{G^I}^{\cal N}(H^I)$, as~claimed.\endproof
 
\thm{firstbig}{Let ${\cal N}$ be a variety of groups, and
let $G\in {\cal N}$ be a finite group. Let $H$ be a subgroup
of~$G$ such that
${\rm dom}_G^{\cal N}(H) = G$. Then
$${\rm dom}_{G\wr\ZZ}^{{\cal N}\!{\cal A}}(H\wr\Z) = G\wr\Z$$
where $\cal A$ is the variety of abelian groups, and $\Z$ is the
infinite cyclic~group.} 

\proof Note that since $H$ is a subgroup of~$G$, $H\wr\Z$ is a
subgroup of~$G\wr\Z$ in the obvious way, by considering only the
functions $\Z\to G$ which take values in~$H$.

By \ref{generalnonsurjinprod}(iii), we need to verify that
$\bigl(H\wr\Z\bigr) {\cal A}(G\wr\Z) = G\wr\Z$, and that
$${\rm dom}_{{\cal A}(G\wr\ZZ)}^{\cal N} \bigl((H\wr\Z)\cap {\cal
A}(G\wr\Z)\bigr)
= {\cal A}(G\wr \Z).$$

We know that for any group $K$, ${\cal A}(K)=[K,K]$, the commutator
of~$K$. By \ref{commofwr}, $[G\wr\Z,G\wr\Z]=G^{\ZZ}$, so $(H\wr\Z){\cal
A}(G\wr\Z)$ contains $G^{\ZZ}$ and also contains $\Z$ (since $\Z$ is a
subgroup of~$H\wr\Z$). As these two subgroups generate $G\wr\Z$, we
conclude that $(H\wr\Z) {\cal A}(G\wr\Z) = G\wr\Z$.

To verify the second condition, note that ${\cal A}(G\wr\Z) =
G^{\ZZ}$, hence we have $$(H\wr\Z)\cap {\cal A}(G\wr\Z) = H^{\ZZ},$$
so we want to find ${\rm dom}_{G^{\ZZ}}^{\cal N}(H^{\ZZ})$. From
\ref{infforfinite} we conclude that $${\rm dom}_{G^{\ZZ}}^{\cal
N}(H^{\ZZ}) = ({\rm dom}_G^{\cal N}(H))^{\ZZ} = G^{\ZZ} = {\cal
A}(G\wr\Z).$$ By \ref{generalnonsurjinprod}(iii), ${\rm dom}_{G\wr\ZZ}^{\cal
NA}(H\wr\Z) = G\wr\Z$, as~claimed.\endproof

\cor{nonsurjtimesab}{Let ${\cal N}$ be a variety of groups
with instances of finite nonsurjective epimorphisms. If ${\cal A}$ is the
variety of abelian groups, then ${\cal NA}$ 
has instances of nonsurjective epimorphisms.\noproof}

\rmrk{caution} Note, however, that the result guarantees that ${\cal
NA}$ has nonsurjective epimorphisms, but does not tell us whether it
also has {\it finite} nonsurjective epimorphisms. Our proof certainly
does not give a finite nonsurjective epimorphism, as $H\wr\Z$ has infinite
index in $G\wr\Z$.

We can extend the result a bit more now:

\lemma{ssubnofwr}{Let $G$ be a group, $n\geq 1$, and ${\cal S}_n$ be the
variety of all solvable groups of solvability length at most~$n$. We
define the groups $K_{G,n}$ recursively, by letting
$$\eqalign{K_{G,1} &= G\wr\Z;\cr
\hbox{and}\qquad K_{G,n+1} &= K_{G,n}\wr\Z.\cr}$$
Then ${\cal S}_n (K_{G,n}) \cong G^{{\ZZ}^n}$.}

\proof We proceed by induction on~$n$. Note that ${\cal
S}_1={\cal A}$, the variety of abelian groups, and ${\cal A}(G\wr\Z)=
[G\wr\Z,G\wr\Z]$; the case $n=1$ now follows from \ref{commofwr}.

Assuming the result is true for $n$, we have
$$\eqalignno{{\cal S}_{n+1}\left(K_{G,n+1}\right) &= {\cal A}\Bigl({\cal
S}_n\bigl(K_{G,n+1}\bigr)\Bigr)\cr
&= {\cal A}\Bigl({\cal S}_n\bigl(K_{G\wr\ZZ,n}\bigr)\Bigr)\cr
&= {\cal A}\Bigl(\bigl(G\wr\Z\bigr)^{{\ZZ}^n}\Bigr)
\qquad\hbox{(by the induction hypothesis)}\cr
&\cong \Bigl({\cal A}(G\wr\Z)\Bigr)^{{\ZZ}^n}
\qquad\hbox{(by \ref{commofwrn})}\cr
&= \bigl(G^{\ZZ}\bigr)^{{\ZZ}^n}\cr
&\cong G^{{\ZZ}^{n+1}}\cr}$$
as~claimed.\endproof

\cor{nonsurjtimesabn}{Let ${\cal N}$ be a variety of groups
with instances of finite nonsurjective epimorphisms, and let $n>0$. Let
${\cal S}_n = {\cal A}^n$ be the variety of all solvable groups of
solvability length at most~$n$. Then ${\cal NS}_n$ has instances
of nonsurjective epimorphisms.}

\proof As before, given any group~$G$, we define recursively
$$\eqalign{ K_{G,1} &= G\wr\Z;\cr
K_{G,n+1} &= K_{G,n}\wr\Z.\cr}$$

Let $G$ be a finite group in ${\cal N}$, and let $H$ a proper subgroup
of~$G$ which is epimorphically embedded into~$G$. We want to show that
$${\rm dom}_{K_{G,n}}^{{\cal S}_n} \bigl(K_{H,n}\bigr)=
K_{G,n}.$$

This follows from \ref{nonsurjinprod}, the fact that ${\cal
S}_n(K_{G,n})\cong G^{{\ZZ}^n}$ (by \ref{ssubnofwr}), and~that $${\rm
dom}_{G^{{\ZZ}^n}}^{\cal N}\left(H^{{\ZZ}^n}\right) = \Bigl({\rm
dom}_{G}^{\cal N}(H)\Bigr)^{{\ZZ}^n}$$
by \ref{infforfinite}.\endproof

Next we prove a result similar to \ref{nonsurjtimesab}
with ${\cal A}$ replaced by an arbitrary variety of nilpotent groups
that contains the infinite cyclic group~$\Z$. 

For $c>0$,  ${\cal N}\!_c$ denotes the variety of nilpotent groups of class at
most~$c$, defined by the single law
$[x_1,x_2,x_3,\cdots,x_{c+1}]$.
In particular, ${\cal N}_1={\cal A}$. We have the following~lemma:

\lemma{verbofnil}{Let $c>0$, and let $G$ be a group. Then
${\cal N}\!_c(G\wr\Z) = G^{\ZZ}$.}

\proof Since every element of $G^{\ZZ}$ may be written in the form
$[\psi,x]$, where $x$ is a generator for $\Z$ and $\psi\in G^{\ZZ}$, by
\ref{commofwr}, it follows that ${\cal
N}_2(G\wr\Z)=G^{\ZZ}$. Proceeding by induction on~$c$ we obtain the
desired~result.\endproof

\thm{secondbig}{Let ${\cal N}$ be a variety of groups,
${\cal N}\!_c$ the variety of nilpotent groups of class at most~$c$,
$c>0$, and let $G\in {\cal N}$ be a finite group. Let $H$ be a
subgroup of~$G$ such that ${\rm dom}_G^{\cal N} (H)= G$. Then
$${\rm dom}_{G\wr\ZZ}^{{\cal N}\!{\cal N}\!_c}(H\wr\Z) = G\wr\Z.$$}

\proof The proof proceeds exactly as the proof of \ref{firstbig},
after we note that
${\cal N}\!_c(G\wr\Z) =G^{\ZZ}$ by \ref{verbofnil}.\endproof

\cor{prodwithnil}{Let ${\cal N}$ be a variety with
instances of finite nonsurjective epimorphisms, and let $c\geq 1$. Then
${\cal N}\!{\cal N}\!_c$ has instances of
nonsurjective~epimorphisms.\noproof}

\cor{prodwithgennil}{Let ${\cal N}$ be a variety with
instances of finite nonsurjective epimorphisms, and let ${\cal Q}$ be a
variety of nilpotent groups such that $\Z\in{\cal Q}$. Then ${\cal
NQ}$ has instances of nonsurjective~epimorphisms.}

\proof Let $G\in {\cal N}$ be a finite group and let $H$ be a
subgroup of~$G$ such that ${\rm dom}_G^{\cal N}(H) = G$. 

Since ${\cal Q}$ is a variety of nilpotent groups, there exists
$c_0\geq 1$ such that ${\cal Q}\subseteq {\cal N}\!_{c_0}$. By
\ref{secondbig}, it follows that ${\rm dom}_{G\wr\ZZ}^{{\cal
NN}\!_{c_{\scriptscriptstyle 0}}}(H\wr\Z)= G\wr\Z$. Since
$G\wr\Z\in{\cal NQ}$, and ${\cal NQ}\subseteq{\cal N}\!{\cal N}\!_{c_0}$,
we~have
$$\eqalign{{\rm dom}_{G\wr\ZZ}^{\cal NQ} (H\wr\Z) &\supseteq {\rm
dom}_{G\wr\ZZ}^{{\cal N}\!{\cal N}\!_{c_{\scriptscriptstyle 0}}}(H\wr\Z)\cr
&= G\wr\Z.\cr}$$
Therefore, ${\cal NQ}$ also has instances of nonsurjective
epimorphisms, as~claimed.\endproof

In the statement of the next Corollary, note that ${\cal N}\!{\cal
N}\!il$ is not, in general, a~variety.

\cor{withacat}{Let ${\cal N}$ be a variety with instances
of finite nonsurjective epimorphisms, and let ${\cal N}\!{\cal N}\!il$ be the
category of all groups which are an extension of an ${\cal N}$-group
by a nilpotent group. Then ${\cal N}\!{\cal N}\!il$ has instances of
nonsurjective epimorphisms.}

\proof Let $G\in{\cal N}$ be a finite group, and $H$ a proper subgroup of~$G$
such that $${\rm dom}_G^{\cal N}(H) = G.$$ We claim that the embedding
$H\wr\Z\hookrightarrow G\wr\Z$ is a nonsurjective epimorphism in the
category~${\cal N}\!{\cal N}\!il$. 

Let $K\in{\cal N}\!{\cal N}\!il$, and
let $f,g\colon (G\wr\Z)\to K$ be two maps which agree on~$H\wr\Z$.

Since $K$ is an extension of an ${\cal N}$-group by a nilpotent group,
there exists a $c>0$ such that $K\in {\cal N}\!{\cal N}\!_c$. But ${\rm
dom}_{G\wr\ZZ}^{{\cal N}\!{\cal N}\!_c}(H\wr\Z)=G\wr\Z$, so
$f|_{H\wr\ZZ}=g|_{H\wr\ZZ}$ implies~$f=g$.

Therefore, the immersion $H\wr\Z\hookrightarrow G\wr\Z$ is an epimorphism in
${\cal N}\!{\cal N}\!il$, which proves the~corollary.\endproof

\thm{bigconclusion}{Let ${\cal N}$ be a variety of groups,
and assume that ${\cal N}$ has instances of finite nonsurjective
epimorphisms. Let ${\cal Q}$ be a variety which is a finite product of
varieties of nilpotent groups, each of which contains the infinite
cyclic group $\Z$. Then ${\cal NQ}$ also has instances of nonsurjective
epimorphisms.}

\proof This is obtained in the same manner as \ref{nonsurjtimesabn}
was obtained before, noting that a formula similar to the one in 
\ref{ssubnofwr} holds for a variety ${\cal Q}$ as~above.\endproof

\rmrk{again} Once again, the condition that ${\cal N}$ have {\it
finite} nonsurjective epimorphisms is somewhat disappointing. It would
be much better if we could drop this extra hypothesis, and I leave as
an open question whether \ref{bigconclusion} holds without~it.

Note that \ref{bigconclusion} implies the weaker
\ref{nonsurjtimesabn}. 

\Section{Another partial answer}{partialtwo}\par

In this section we will prove another partial converse to
\ref{nsinprodonlyiffirst}. Namely, we will show that if a variety
${\cal N}$ contains a nonsurjective epimorphism of the form $H\to S$,
where $S$ is a finite nonabelian simple group, then for any nontrivial
variety~${\cal Q}$, the variety ${\cal NQ}$ also contains a (finite)
nonsurjective epimorphism.

The idea behind the argument is simple. Suppose that we have a variety
${\cal N}$, a finite simple nonabelian group~$S$ in~${\cal N}$, and
a proper subgroup $H$ of~$S$ which is epimorphically embedded into~$S$
in the variety~${\cal N}$.

If ${\cal Q}$ is any variety, then either $S\in {\cal Q}$ or $S\notin
{\cal Q}$. If $S\notin{\cal Q}$, then ${\cal Q}(S)=S$ (since
$S$ is simple) and so by \ref{nonsurjinprod} the embedding $H\to S$ is
also a nonsurjective epimorphism in~${\cal NQ}$. On the other hand, if
$S\in {\cal Q}$, we want to find some finite group $G\in {\cal Q}$
such that $S\wr G$ is not in ${\cal Q}$. We might expect such
a~$G$ to exist, say some finite group which is ``barely''
in~${\cal Q}$. Then we can hope that 
$H\wr G$ will be epimorphically embedded into $S\wr G$ in~${\cal
NQ}$.

Indeed, it turns out that this is the case. We need a few preparatory
lemmas to establish the existence of a~$G$ with the properties above.

\lemma{notallpgroups}{Let ${\cal V}$ be a variety of groups, and $p$
a prime. If ${\cal V}$ contains all finite $p$-groups, then ${\cal V}$
is the variety of all groups.}

\proof First, note that an absolutely free group can be embedded in a
noncommuting formal power series ring, as the formal power series with
constant term $1$ over $\Z/p\Z$. Namely, the absolutely free group on
$x_1,\ldots,x_n$ can be embedded into $\Z_p\langle\!\langle
y_1,\ldots,y_n\rangle\!\rangle$ by sending $x_i$ into $1+y_i$. One
shows this is an embedding by noting that if $n>0$ is an integer, and
$n=ap^m$, where ${\rm gcd}(a,p)=1$,~then
$$(1+y_i)^m = 1 + ay_i^{p^m} + y_i^{2p^m}h(y_i)$$
where $h(y_i)$ is a power series on~$y_i$. Then, if
${\bf w}=x_{r_1}^{a_1}x_{r_2}^{a_2}\cdots x_{r_s}^{a_s}$ is a
nontrivial reduced word on
the $x_i$'s (so that $r_i\not= r_{i+1}$, and $a_j\not=0$ for all $i$
and~$j$), if we write $A_i = b_ip^{k_i}$, with ${\rm gcd}(b_i,p)=1$, then
the image of this word in $\Z_p\langle\!\langle
y_1,\ldots,y_n\rangle\!\rangle$ under the map given above has a unique
monomial of degree $p^{k_1}+\cdots+p^{k_s}$, namely
$$a_1\cdots a_sy_{r_1}^{p^{k_1}}\cdots y_{r_s}^{p^{k_s}}$$ (see 
Theorem~5.6 in~{\bf [\cite{magnus}]}, for the~details).

Next, note that every nontrivial identity on $n$ variables will fail in a
truncated formal power series ring in finitely many noncommuting
indeterminates, 
$$\Z_p\langle\!\langle y_1,\ldots,y_n\rangle\!\rangle/(y_1,\ldots,y_n)^d$$
for some $d$, (namely, if the word is ${\bf w}$ as above, setting
$d>p^{k_1} + \cdots + p^{k_n}$ will guarantee that the image of ${\bf
w}$ is nontrivial). And the set of such truncated power series with
constant term 1 is a finite $p$ group, since it has $p^{d-1}$ elements.

Therefore, if ${\cal V}$ contains all finite $p$ groups, then no
nontrivial identity can be a law of~${\cal V}$, from which it follows
that ${\cal V}$ is the variety of all groups.\endproof
 
\cor{icantkeepgoing}{Let ${\cal V}$ be a proper subvariety of
${\cal G}roup$. Then there exists a finite group $G\in {\cal V}$ such that
$(\Z/p\Z)\wr G\notin {\cal V}$.}

\proof Suppose that for all finite groups $G\in {\cal V}$,
$(\Z/p\Z)\wr G$ is also in ${\cal V}$. 

Every finite $p$ group can be obtained from the trivial subgroup by
successively extending $\Z/p\Z$. To see this, recall that a finite $p$
group always has nontrivial center, so the group must contain a normal
subgroup isomorphic to $\Z/p\Z$. Since any extension of $\Z/p\Z$ by
$G$ can be realized as a subgroup of $(\Z/p\Z)\wr G$ by a theorem of
Kaloujnine and Krasner {\bf [\cite{wreathext}]}, it follows that
${\cal V}$ must contain all finite $p$ groups. But in that case,
\ref{notallpgroups} shows that ${\cal V}$ is the variety of all
groups.\endproof

\cor{finiteinvarworks}{Let ${\cal V}$ be a nontrivial variety of
groups, and let $A\in {\cal V}$ be a nontrivial group. Then there
exists a finite group~$G\in {\cal V}$ such that $A\wr G\notin{\cal V}$.}

\proof Let $x$ be a nontrivial element of~$A$, and let $p$ be a prime
such that $\langle x\rangle$ has $\Z/p\Z$ as a homomorphic image. Then
for every group~$B$, $(\Z/p\Z)\wr B$ is a homomorphic image of a
subgroup of~$A\wr
B$. By \ref{icantkeepgoing} there exists a finite
group $G$ such that $(\Z/p\Z)\wr G\notin{\cal V}$. In particular $A\wr
G$ cannot be in~${\cal V}$.\endproof

The final ingredient for the proof is the following lemma:

\lemma{qofsimplewrB}{Let ${\cal N}$ be a variety, and let $S$ be a
finite nonabelian simple group in ${\cal N}$. Let ${\cal Q}$ be a
nontrivial variety, and let $B$ be a finite group in ${\cal Q}$. Then
either ${\cal Q}(S\wr B) = S^B$ or~${\cal Q}(S\wr B) = \{e\}$.}

\proof Since $(S\wr B)/S^B\cong B\in {\cal Q}$, we have ${\cal Q}(S\wr
B)\subseteq S^B.$ Also ${\cal Q}(S\wr B)\triangleleft S\wr B$, and so
in particular ${\cal Q}(S\wr B)\triangleleft S^B$. Since $S$ is a
finite nonabelian simple group, and $B$ is a finite group, 
${\cal Q}(S\wr B)$ must equal the product of
some of the copies of~$S$ (see for example {\bf [\cite{remak}]}).

Now we simply note that in $S\wr B$, $B$ acts transitively on
the factors $S$ of $S^B$, hence given a subgroup of~$S^B$ which is the
product of a subset of these factors, and is invariant under the
action of~$B$ (as ${\cal Q}(S\wr B)$ must be), it must be the product
of all the factors, or of none.\endproof

\thm{simpletimesanything}{Let ${\cal N}$ be a nontrivial variety, $S$
a finite nonabelian simple group, and $H$ a proper subgroup of~$S$
such that $S\in {\cal N}$ and ${\rm dom}_S^{{\cal N}}(H)=S$ (that
is, the embedding $H\hookrightarrow S$ is a nonsurjective epimorphism
in ${\cal N}$). Let ${\cal Q}$ be any nontrivial variety. Then ${\cal
NQ}$ also has instances of finite nonsurjective epimorphisms. Namely,
there exists a finite (possibly trivial) group $G\in {\cal Q}$ such
that $H\wr G$ is epimorphically embedded in $S\wr G$, in the
variety~${\cal NQ}$.}

\proof Let $H$, $S$, ${\cal N}$ and ${\cal Q}$ be as in the statement.
By \ref{finiteinvarworks} there exists a finite group $G \in {\cal Q}$
such that $S\wr G\notin{\cal Q}$ ($G$ may be trivial, for example if
$S\notin {\cal Q}$). By \ref{qofsimplewrB}, we have~${\cal Q}(S\wr G)
= S^G$, since $S\wr G\notin {\cal Q}$.

Consider the subgroup $H\wr G$ of $S\wr G$. $(H\wr G)\,{\cal Q}(S\wr
G)$ equals the whole group, 
since it contains both $S^G$ and~$G$, which together generate $S\wr G$.
On the other hand, 
$$\bigl(H\wr G\bigr) \cap {\cal Q}(S\wr G) = \bigl(H\wr G\bigr) \cap S^G
= H^G.$$
Since dominions respect finite direct products,
$$\eqalign{{\rm dom}_{S^G}^{\cal N}\Bigl(H^G\Bigr) &= \Bigl({\rm
dom}_S^{\cal N}(H)\Bigr)^G\cr &= S^G\cr &= {\cal Q}(S\wr G).\cr}$$ By
\ref{generalnonsurjinprod}(iii), ${\rm dom}_{S\wr G}^{\cal NQ}(H\wr G)
= S\wr G$, so the embedding $H\wr G \hookrightarrow S\wr G$ is a
nonsurjective epimorphism in ${\cal NQ}$. Since all the groups
involved are finite, it is an instance of a finite nonsurjective
epimorphism, as~claimed.\endproof

This provides a nice partial answer to the question at the end of
\ref{charactofepis}, especially since all basic examples we have given
so far (for example, in {\bf [\cite{simpleprelim}]}) 
are precisely of the
kind described in \ref{simpletimesanything}, namely the embedding of a
proper subgroup into a simple nonabelian group. Since we also know
that in a variety consisting only of solvable groups all epimorphisms
are surjective, we might guess that finite nonabelian simple groups
will play a large role whenever a nonsurjective epimorphism occurs.

We finish this section with the example promised in
\ref{needallnotjustone}.

{\bf Example~\numbeq{afivewrafivesquared}.} Let ${\cal N}={\cal
Q}={\rm Var}(A_5\wr A_5)$, and let $G=(A_5\wr A_5)\wr A_5\in {\cal
NQ}$. Let $N=(A_5\wr A_5)^{A_5}$ (that is, the base group of the last
wreath product taken in the construction of~$G$), and let $H$ be the
subgroup of~$G$ given by $H=(A_4\wr A_5)\wr A_5$. Then $N\in {\cal
N}$, $G/N\cong A_5\in {\cal Q}$, $NH=G$, and $H\cap N=(A_4\wr
A_5)^{A_5}$, which is epimorphically embedded into $N$ in ${\cal
N}$. However, $H$ is not epimorphically embedded in~$G$ in the variety
${\cal NQ}$. To see this, we use the characterization in
\ref{nonsurjinprod}. Note that ${\cal Q}(G)$ is the subgroup
$(A_5^{A_5})^{A_5}$, so $H\cap {\cal B}(G)$ is $(A_4^{A_5})^{A_5}$; we
know $A_4$ is not epimorphically embedded into~$A_5$ in ${\cal N}$, so
$H\cap {\cal Q}(G)$ is not epimorphically embedded into ${\cal
Q}(G)$. In particular, $H\hookrightarrow G$ is not an epimorphism in~${\cal
NQ}$. This shows that in \ref{generalnonsurjinprod}(iv), it is not
enough to consider a single normal subgroup~$N$.

%
\ifnum0<\citations{\par\bigbreak
\filbreak{\bf References}\par\frenchspacing}\fi
%
\ifundefined{xthreeNB}\else
\item{\bf [\refer{threeNB}]}{Baumslag, G{.}, Neumann, B{.}H{.},
Neumann, H{.}, and Neumann, P{.}M. {\it On varieties generated by a
finitely generated group.\/} {\sl Math.\ Z.} {\bf 86} (1964)
pp.~\hbox{93--122}. {MR:30\#138}}\par\filbreak\fi
\ifundefined{xbergman}\else
\item{\bf [\refer{bergman}]}{Bergman, George M. {\it An Invitation to
General Algebra and Universal Constructions.\/} {\sl Berkeley
Mathematics Lecture Notes 7\/} (1995).}\par\filbreak\fi
\ifundefined{xordersberg}\else
\item{\bf [\refer{ordersberg}]}{Bergman, George M. {\it Ordering
coproducts of groups and semigroups.\/} {\sl J. Algebra} {\bf 133} (1990)
no. 2, pp.~\hbox{313--339}. {MR:91j:06035}}\par\filbreak\fi
\ifundefined{xbirkhoff}\else
\item{\bf [\refer{birkhoff}]}{Birkhoff, Garrett. {\it On the structure
of abstract algebras.\/} {\sl Proc.\ Cambridge\ Philos.\ Soc.} {\bf
31} (1935), pp.~\hbox{433--454}.}\par\filbreak\fi
\ifundefined{xbrown}\else
\item{\bf [\refer{brown}]}{Brown, Kenneth S. {\it Cohomology of
Groups, 2nd Edition.\/} {\sl Graduate texts in mathematics 87\/},
Springer Verlag,~1994. {MR:96a:20072}}\par\filbreak\fi
\ifundefined{xmetab}\else
\item{\bf [\refer{metab}]}{Golovin, O. N. {\it Metabelian products of
groups.\/}
{\sl American Mathematical Society Translations}, series 2, {\bf 2} (1956),
pp.~\hbox{117--131.} {MR:17,824b}}\par\filbreak\fi
\ifundefined{xhall}\else
\item{\bf [\refer{hall}]}{Hall, M. {\it The Theory of Groups.\/}
Mac~Millan Company,~1959. {MR:21\#1996}}\par\filbreak\fi
\ifundefined{xphall}\else
\item{\bf [\refer{phall}]}{Hall, P. {\it Verbal and marginal
subgroups.} {\sl J.\ Reine\ Angew.\ Math.\/} {\bf 182} (1940)
pp.~\hbox{156--157.} {MR:2,125i}}\par\filbreak\fi
\ifundefined{xheineken}\else
\item{\bf [\refer{heineken}]}{Heineken, H. {\it Engelsche Elemente der
L\"ange drei,\/} {\sl Illinois Journal of Math.} {\bf 5} (1961)
pp.~\hbox{681--707.} {MR:24\#A1319}}\par\filbreak\fi
\ifundefined{xherman}\else
\item{\bf [\refer{herman}]}{Herman, Krzysztof. {\it Some remarks on
the twelfth problem of Hanna Neumann.\/} {\sl Publ.\ Math.\ Debrecen}
{\bf 37} (1990)  no. 1--2, pp.~\hbox{25--31.} {MR:91f:20030}}\par\filbreak\fi
\ifundefined{xherstein}\else
\item{\bf [\refer{herstein}]}{Herstein, I.~N. {\it Topics in
Algebra.\/} Blaisdell Publishing Co.,~1964.}\par\filbreak\fi
\ifundefined{xepisandamalgs}\else
\item{\bf [\refer{episandamalgs}]}{Higgins, Peter M. {\it Epimorphisms
and amalgams.} {\sl
Colloq.\ Math.} {\bf 56} no.~1 (1988) pp.~\hbox{1--17.}
{MR:89m:20083}}\par\filbreak\fi
\ifundefined{xhigmanpgroups}\else
\item{\bf [\refer{higmanpgroups}]}{Higman, Graham. {\it Amalgams of
$p$-groups.\/} {\sl J. of~Algebra} {\bf 1} (1964)
pp.~\hbox{301--305.} {MR:29\#4799}}\par\filbreak\fi
\ifundefined{xhigmanremarks}\else
\item{\bf [\refer{higmanremarks}]}{Higman, Graham. {\it Some remarks
on varieties of groups.\/} {\sl Quart.\ J.\ of Math.\ (Oxford) (2)} {\bf
10} (1959), pp.~\hbox{165--178.} {MR:22\#4756}}\par\filbreak\fi
\ifundefined{xhughes}\else
\item{\bf [\refer{hughes}]}{Hughes, N.J.S. {\it The use of bilinear
mappings in the classification of groups of class~$2$.\/} {\sl Proc.\
Amer.\ Math.\ Soc.\ } {\bf 2} (1951) pp.~\hbox{742--747.}
{MR:13,528e}}\par\filbreak\fi
\ifundefined{xisbelltwo}\else
\item{\bf [\refer{isbelltwo}]}{Howie, J.~M., Isbell, J.~R. {\it
Epimorphisms and dominions II.\/} {\sl Journal of Algebra {\bf
6}}(1967) pp.~\hbox{7--21.} {MR:35\#105b}}\par\filbreak\fi
\ifundefined{xisaacs}\else
\item{\bf [\refer{isaacs}]}{Isaacs, I.M., Navarro, Gabriel. {\it
Coprime actions, fixed-point subgroups and irreducible induced
characters.} {\sl J.~of Algebra} {\bf 185} (1996) no.~1,
pp.~\hbox{125--143.} {MR:97g:20009}}\par\filbreak\fi
\ifundefined{xisbellone}\else
\item{\bf [\refer{isbellone}]}{Isbell, J. R. {\it Epimorphisms and
dominions} in {\sl 
Proc.~of the Conference on Categorical Algebra, La Jolla 1965,\/}
pp.~\hbox{232--246.} Lange and Springer, New
York~1966. MR:35\#105a (The statement of the
Zigzag Lemma for {\it rings} in this paper is incorrect. The correct
version is stated in~{\bf [\cite{isbellfour}]})}\par\filbreak\fi
\ifundefined{xisbellthree}\else
\item{\bf [\refer{isbellthree}]}{Isbell, J. R. {\it Epimorphisms and
dominions III.} {\sl Amer.\ J.\ Math.\ }{\bf 90} (1968)
pp.~\hbox{1025--1030.} {MR:38\#5877}}\par\filbreak\fi
\ifundefined{xisbellfour}\else
\item{\bf [\refer{isbellfour}]}{Isbell, J. R. {\it Epimorphisms and
dominions IV.} {\sl Journal\ London Math.\ Society~(2),}
{\bf 1} (1969) pp.~\hbox{265--273.} {MR:41\#1774}}\par\filbreak\fi
\ifundefined{xjones}\else
\item{\bf [\refer{jones}]}{Jones, Gareth A. {\it Varieties and simple
groups.\/} {\sl J.\ Austral.\ Math.\ Soc.} {\bf 17} (1974)
pp.~\hbox{163--173.} {MR:49\#9081}}\par\filbreak\fi
\ifundefined{xjonsson}\else
\item{\bf [\refer{jonsson}]}{J\'onsson, B. {\it Varieties of groups of
nilpotency three.} {\sl Notices Amer.\ Math.\ Soc.} {\bf 13} (1966)
pp.~488.}\par\filbreak\fi
\ifundefined{xwreathext}\else
\item{\bf [\refer{wreathext}]}{Kaloujnine, L. and Krasner, Marc. {\it
Produit complet des groupes de permutations et le probl\`eme
d'extension des groupes III.} {\sl Acta Sci.\ Math.\ Szeged} {\bf 14}
(1951) pp.~\hbox{69--82}. {MR:14,242d}}\par\filbreak\fi
\ifundefined{xkhukhro}\else
\item{\bf [\refer{khukhro}]}{Khukhro, Evgenii I. {\it Nilpotent Groups
and their Automorphisms.} {\sl de Gruyter Expositions in Mathematics}
{\bf 8}, New York 1993. {MR:94g:20046}}\par\filbreak\fi
\ifundefined{xkleimanbig}\else
\item{\bf [\refer{kleimanbig}]}{Kle\u{\i}man, Yu.~G. {\it On
identities in groups.\/} {\sl Trans.\ Moscow Math.\ Soc.\ } 1983,
Issue 2, pp.~\hbox{63--110}. {MR:84e:20040}}\par\filbreak\fi
\ifundefined{xthirtynine}\else
\item{\bf [\refer{thirtynine}]}{Kov\'acs, L.~G. {\it The thirty-nine
varieties.} {\sl Math.\ Scientist} {\bf 4} (1979)
pp.~\hbox{113--128.} {MR:81m:20037}}\par\filbreak\fi
\ifundefined{xlamssix}\else
\item{\bf [\refer{lamssix}]}{Lam, T{.}Y{.}, and Leep, David B. {\it
Combinatorial structure on the automorphism group of~$S_6$.\/} {\sl
Expo. Math.} {\bf 11} (1993) pp.~\hbox{289--308.}
{MR:94i:20006}}\par\filbreak\fi
\ifundefined{xlevione}\else
\item{\bf [\refer{levione}]}{Levi, F.~W. {\it Groups on which the
commutator relation 
satisfies certain algebraic conditions.\/} {\sl J.\ Indian Math.\ Soc.\ New
Series} {\bf 6}(1942), pp.~\hbox{87--97.} {MR:4,133i}}\par\filbreak\fi
\ifundefined{xgermanlevi}\else
\item{\bf [\refer{germanlevi}]}{Levi, F.~W. and van der Waerden,
B.~L. {\it \"Uber eine 
besondere Klasse von Gruppen.\/} {\sl Abhandl.\ Math.\ Sem.\ Univ.\ Hamburg}
{\bf 9}(1932), pp.~\hbox{154--158.}}\par\filbreak\fi
\ifundefined{xlichtman}\else
\item{\bf [\refer{lichtman}]}{Lichtman, A.~L. {\it Necessary and
sufficient conditions for the residual nilpotence of free products of
groups.\/} {\sl J. Pure and Applied Algebra} {\bf 12} no. 1 (1978),
pp.~\hbox{49--64.} {MR:58\#5938}}\par\filbreak\fi
\ifundefined{xmaxofan}\else
\item{\bf [\refer{maxofan}]}{Liebeck, Martin W.; Praeger, Cheryl E.;
and Saxl, Jan. {\it A classification of the maximal subgroups of the
finite alternating and symmetric groups.\/} {\sl J. of Algebra} {\bf
111}(1987), pp.~\hbox{365--383.} {MR:89b:20008}}\par\filbreak\fi
\ifundefined{xepisingroups}\else
\item{\bf [\refer{episingroups}]}{Linderholm, C.E. {\it A group
epimorphism is surjective.\/} {\sl Amer.\ Math.\ Monthly\ }77
pp.~\hbox{176--177.}}\par\filbreak\fi
\ifundefined{xmckay}\else
\item{\bf [\refer{mckay}]}{McKay, Susan. {\it Surjective epimorphisms
in classes
of groups.} {\sl Quart.\ J.\ Math.\ Oxford (2),\/} {\bf 20} (1969),
pp.~\hbox{87--90.} {MR:39\#1558}}\par\filbreak\fi
\ifundefined{xmaclane}\else
\item{\bf [\refer{maclane}]}{Mac Lane, Saunders. {\it Categories for
the Working Mathematician.} {\sl Graduate texts in mathematics 5},
Springer Verlag (1971). {MR:50\#7275}}\par\filbreak\fi
\ifundefined{xbilinear}\else
\item{\bf [\refer{bilinear}]}{Magidin, Arturo. {\it Bilinear maps and 
2-nilpotent groups.\/} August 1996, 7~pp.}\par\filbreak\fi
\ifundefined{xbilinearprelim}\else
\item{\bf [\refer{bilinearprelim}]}{Magidin, Arturo. {\it Bilinear maps
and central extensions of abelian groups.\/} In~preparation.}\par\filbreak\fi
\ifundefined{xprodvar}\else
\item{\bf [\refer{prodvar}]}{Magidin, Arturo. {\it Dominions in product
varieties of groups.\/} May 1997, 21~pp.}\par\filbreak\fi
\ifundefined{xprodvarprelim}\else
\item{\bf [\refer{prodvarprelim}]}{Magidin, Arturo. {\it Dominions in product
varieties of groups.\/} In preparation.}\par\filbreak\fi
\ifundefined{xmythesis}\else
\item{\bf [\refer{mythesis}]}{Magidin, Arturo. {\it Dominions in
Varieties of Groups.\/} Doctoral dissertation, University of
California at Berkeley, May 1998.}\par\filbreak\fi
\ifundefined{xnildoms}\else
\item{\bf [\refer{nildoms}]}{Magidin, Arturo {\it Dominions in varieties
of nilpotent groups.\/} December 1996, 27~pp.}\par\filbreak\fi
\ifundefined{xnildomsprelim}\else
\item{\bf [\refer{nildomsprelim}]}{Magidin, Arturo. {\it Dominions in
varieties of nilpotent groups.\/} In preparation.}\par\filbreak\fi
\ifundefined{xsimpleprelim}\else
\item{\bf [\refer{simpleprelim}]}{Magidin, Arturo. {\it Dominions in
varieties generated by simple groups.\/} In preparation.}\par\filbreak\fi
\ifundefined{xntwodoms}\else
\item{\bf [\refer{ntwodoms}]}{Magidin, Arturo. {\it Dominions in the variety of
2-nilpotent groups.\/} May 1996, 6~pp.}\par\filbreak\fi
\ifundefined{xdomsmetabprelim}\else
\item{\bf [\refer{domsmetabprelim}]}{Magidin, Arturo. {\it Dominions
in the variety of metabelian groups.\/}
In~preparation.}\par\filbreak\fi
\ifundefined{xfgnilgroups}\else
\item{\bf [\refer{fgnilgroups}]}{Magidin, Arturo. {\it Dominions of
finitely generated nilpotent groups.\/} October~1997,
10~pp.}\par\filbreak\fi
\ifundefined{xfgnilprelim}\else
\item{\bf [\refer{fgnilprelim}]}{Magidin, Arturo. {\it Dominions of
finitely generated nilpotent groups.\/} In preparation.}\par\filbreak\fi
\ifundefined{xwordsprelim}\else
\item{\bf [\refer{wordsprelim}]}{Magidin, Arturo. {\it
Words and dominions.\/} In~preparation.}\par\filbreak\fi
\ifundefined{xepis}\else
\item{\bf [\refer{epis}]}{Magidin, Arturo. {\it Non-surjective epimorphisms
in varieties of groups and other results.\/} February 1997,
13~pp.}\par\filbreak\fi
\ifundefined{xoddsandends}\else
\item{\bf [\refer{oddsandends}]}{Magidin, Arturo. {\it Some odds and
ends.\/} June 1996, 3~pp.}\par\filbreak\fi
\ifundefined{xpropdom}\else
\item{\bf [\refer{propdom}]}{Magidin, Arturo. {\it Some properties of
dominions in varieties of groups.\/} March 1997, 13~pp.}\par\filbreak\fi
\ifundefined{xzabsp}\else
\item{\bf [\refer{zabsp}]}{Magidin, Arturo. {\it $\Z$ is an absolutely
closed $2$-nil group.\/} Submitted.}\par\filbreak\fi
\ifundefined{xmagnus}\else
\item{\bf [\refer{magnus}]}{Magnus, Wilhelm; Karras, Abraham; and
Solitar, Donald. {\it Combinatorial Group Theory.\/} 2nd Edition; Dover
Publications, Inc.~1976. {MR:53\#10423}}\par\filbreak\fi
\ifundefined{xamalgtwo}\else
\item{\bf [\refer{amalgtwo}]}{Maier, Berthold J. {\it Amalgame
nilpotenter Gruppen
der Klasse zwei II.\/} {\sl Publ.\ Math.\ Debrecen} {\bf 33}(1986),
pp.~\hbox{43--52.} {MR:87k:20050}}\par\filbreak\fi
\ifundefined{xnilexpp}\else
\item{\bf [\refer{nilexpp}]}{Maier, Berthold J. {\it On nilpotent
groups of exponent $p$.\/} {\sl Journal of~Algebra} {\bf 127} (1989)
pp.~\hbox{279--289.} {MR:91b:20046}}\par\filbreak\fi
\ifundefined{xmaltsev}\else
\item{\bf [\refer{maltsev}]}{Maltsev, A.~I. {\it Generalized
nilpotent algebras and their associated groups.} (Russian) {\sl
Mat.\ Sbornik N.S.} {\bf 25(67)} (1949) pp.~\hbox{347--366.} ({\sl
Amer.\ Math.\ Soc.\ Translations Series 2} {\bf 69} 1968,
pp.~\hbox{1--21.}) {MR:11,323b}}\par\filbreak\fi
\ifundefined{xmaltsevtwo}\else
\item{\bf [\refer{maltsevtwo}]}{Maltsev, A.~I. {\it Homomorphisms onto
finite groups.} (Russian) {\sl Ivanov. gosudarst. ped. Inst., u\v
cenye zap., fiz-mat. Nuak} {\bf 18} (1958)
\hbox{pp. 49--60.}}\par\filbreak\fi
\ifundefined{xmorandual}\else
\item{\bf [\refer{morandual}]}{Moran, S. {\it Duals of a verbal
subgroup.\/} {\sl J.\ London Math.\ Soc.} {\bf 33} (1958)
pp.~\hbox{220--236.} {MR:20\#3909}}\par\filbreak\fi
\ifundefined{xhneumann}\else
\item{\bf [\refer{hneumann}]}{Neumann, Hanna. {\it Varieties of
Groups.\/} {\sl Ergebnisse der Mathematik und ihrer Grenz\-ge\-biete\/}
New series, Vol.~37, Springer Verlag~1967. {MR:35\#6734}}\par\filbreak\fi
\ifundefined{xneumannwreath}\else
\item{\bf [\refer{neumannwreath}]}{Neumann, Peter M. {\it On the
structure of standard wreath products of groups.\/} {\sl Math.\
Zeitschr.\ }{\bf 84} (1964) pp.~\hbox{343--373.} {MR:32\#5719}}\par\filbreak\fi
\ifundefined{xpneumann}\else
\item{\bf [\refer{pneumann}]}{Neumann, Peter M. {\it Splitting groups
and projectives
in varieties of groups.\/} {\sl Quart.\ J.\ Math.\ Oxford} (2), {\bf
18} (1967),
pp.~\hbox{325--332.} {MR:36\#3859}}\par\filbreak\fi
\ifundefined{xoates}\else
\item{\bf [\refer{oates}]}{Oates, Sheila. {\it Identical Relations in
Groups.\/} {\sl J.\ London Math.\ Soc.} {\bf 38} (1963),
pp.~\hbox{71--78.} {MR:26\#5043}}\par\filbreak\fi
\ifundefined{xolsanskii}\else
\item{\bf [\refer{olsanskii}]}{Ol'\v{s}anski\v{\i}, A. Ju. {\it On the
problem of a finite basis of identities in groups.\/} {\sl
Izv.\ Akad.\ Nauk.\ SSSR} {\bf 4} (1970) no. 2
pp.~\hbox{381--389.}}\par\filbreak\fi
\ifundefined{xremak}\else
\item{\bf [\refer{remak}]}{Remak, R. {\it \"Uber minimale invariante
Untergruppen in der Theorie der end\-lichen Gruppen.\/} {\sl
J.\ reine.\ angew.\ Math.} {\bf 162} (1930),
pp.~\hbox{1--16.}}\par\filbreak\fi
\ifundefined{xclassifthree}\else
\item{\bf [\refer{classifthree}]}{Remeslennikov, V. N. {\it Two
remarks on 3-step nilpotent groups} (Russian) {\sl Algebra i Logika
Sem.} (1965) no.~2 pp.~\hbox{59--65.} {MR:31\#4838}}\par\filbreak\fi
\ifundefined{xrotman}\else
\item{\bf [\refer{rotman}]}{Rotman, J.J. {\it Introduction to the Theory of
Groups}, 4th edition. {\sl Graduate texts in mathematics 119},
Springer Verlag,~1994. {MR:95m:20001}}\par\filbreak\fi
\ifundefined{xsaracino}\else
\item{\bf [\refer{saracino}]}{Saracino, D. {\it Amalgamation bases for
nil-$2$ groups.\/} {\sl Alg.\ Universalis\/} {\bf 16} (1983),
pp.~\hbox{47--62.} {MR:84i:20035}}\par\filbreak\fi
\ifundefined{xscott}\else
\item{\bf [\refer{scott}]}{Scott, W.R. {\it Group Theory.} Prentice
Hall,~1964. {MR:29\#4785}}\par\filbreak\fi
\ifundefined{xsmelkin}\else
\item{\bf [\refer{smelkin}]}{\v{S}mel'kin, A.L. {\it Wreath products and
varieties of groups} [Russian] {\sl Dokl.\ Akad.\ Nauk S.S.S.R.\/} {\bf
157} (1964), pp.~\hbox{1063--1065} Transl.: {\sl Soviet Math.\ Dokl.\ } {\bf
5} (1964), pp.~\hbox{1099--1011}. {MR:33\#1352}}\par\filbreak\fi
\ifundefined{xstruikone}\else
\item{\bf [\refer{struikone}]}{Struik, Ruth Rebekka. {\it On nilpotent
products of cyclic groups.\/} {\sl Canadian Journal of
Mathematics\/} {\bf 12} (1960)
pp.~\hbox{447--462}. {MR:22\#11028}}\par\filbreak\fi
\ifundefined{xstruiktwo}\else
\item{\bf [\refer{struiktwo}]}{Struik, Ruth Rebekka. {\it On nilpotent
products of cyclic groups II.\/} {\sl Canadian Journal of
Mathematics\/} {\bf 13} (1961) pp.~\hbox{557--568.}
{MR:26\#2486}}\par\filbreak\fi
\ifundefined{xvlee}\else
\item{\bf [\refer{vlee}]}{Vaughan-Lee, M{.} R{.} {\it Uncountably many
varieties of groups.\/} {\sl Bull.\ London Math.\ Soc.} {\bf 2} (1970)
pp.~\hbox{280--286.} {MR:43\#2054}}\par\filbreak\fi
\ifundefined{xweibel}\else
\item{\bf [\refer{weibel}]}{Weibel, Charles. {\it Introduction to
Homological Algebra.\/} Cambridge University
Press~1994. {MR:95f:18001}}\par\filbreak\fi 
\ifundefined{xweigelone}\else
\item{\bf [\refer{weigelone}]}{Weigel, T.S. {\it Residual properties
of free groups.\/} {\sl J.\ of Algebra} {\bf 160} (1993)
pp.~\hbox{14--41.} {MR:94f:20058a}}\par\filbreak\fi
\ifundefined{xweigeltwo}\else
\item{\bf [\refer{weigeltwo}]}{Weigel, T.S. {\it Residual properties
of free groups II.\/} {\sl Comm.\ in Algebra} {\bf 20}(5) (1992)
pp.~\hbox{1395--1425.} {MR:94f:20058b}}\par\filbreak\fi
\ifundefined{xweigelthree}\else 
\item{\bf [\refer{weigelthree}]}{Weigel, T.S. {\it Residual Properties
of free groups III.\/} {\sl Israel J.\ Math.\ } {\bf 77} (1992)
pp.~\hbox{65--81.} {MR:94f:20058c}}\par\filbreak\fi
\ifundefined{xzstwo}\else
\item{\bf [\refer{zstwo}]}{Zariski, Oscar and Samuel, Pierre. {\it
Commutative Algebra}, Volume
II. Springer-Verlag~1976. {MR:52\#10706}}\par\filbreak\fi
\ifnum0<\citations\nonfrenchspacing\fi

\bigskip
{\it
\obeylines
\noindent Arturo Magidin
\noindent Cub\'iculo 112
\noindent Instituto de Matem\'aticas
\noindent Universidad Nacional Aut\'onoma de M\'exico
\noindent 04510 Mexico City, MEXICO
\noindent e-mail: magidin@matem.unam.mx
}
 
\vfill\eject
\immediate\closeout\aux
\end